\documentclass[10pt,twoside]{classe-ma_f}

\usepackage{amssymb,amsbsy,amsmath,amsfonts,epsfig}
\usepackage[english,francais]{babel}
\usepackage{times}

 \def\ZZ{{\mathbb{Z}}}
\def\og{\leavevmode\raise.3ex\hbox{$\scriptscriptstyle\langle\!\langle$~}}
\def\fg{\leavevmode\raise.3ex\hbox{~$\!\scriptscriptstyle\,\rangle\!\rangle$}}
\DeclareMathOperator{\Hom}{Hom}     %homomorphisms
\newcommand{\mm}{\underline{m}}
\newtheorem{e-proposition}[theorem]{Proposition}

\newtheorem{e-definition}[theorem]{Definition\rm}

\newtheorem{theoreme}{Th\'eor\`eme}[section]

\newtheorem{corollaire}[theoreme]{Corollaire}

\newtheorem{exemple}{\it Exemple}

\ComParit{S0764-4442}
\PIT{FLA}
\PXMA{2159}
\Add{0}
\Volume{333}
\Year{2001}
\FirstPage{939}
\LastPage{942}
\AuteurCourant{David Cimasoni}
\TitreCourant{L'homologie de Novikov des entrelacs de Waldhausen}
\Journal
\Rubrique{Topologie}{Topology}
\PresentePar{\'Etienne}{Ghys}
\Recu{2 octobre 2001}{}{8 octobre 2001}
\begin{document}
\selectlanguage{francais}
\setcounter{page}{939}
\title{%
L'homologie de Novikov des entrelacs de Waldhausen
}
\author{%
David Cimasoni
}
\address{%
Section de math\'ematiques, Universit\'e de Gen\`eve, 2-4 rue du Li\`evre, 1211 Gen\`eve 24, Suisse\\
Courriel:David.Cimasoni@math.unige.ch
}
\maketitle
\thispagestyle{empty}
%%%%%%%%%%%%%%%%%%%%%%%%%%%%%%%%%%%%%%%%%%%%%%%%%%%%%%%%%%%%
%%%  R\'esum\'e  %%%
%%%%%%%%%%%%%%%%%%%%
\begin{Resume}{%
Un multi-entrelacs de Waldhausen est un entrelacs avec multiplicit\'es dans une 3-sph\`ere
d'homologie dont l'ext\'erieur est une vari\'et\'e de Waldhausen. Dans cette Note, nous
calculons l'homologie de Novikov des multi-entrelacs de Waldhausen.
En corollaire, nous obtenons une borne pour le nombre de modules de Novikov d'un entrelacs de
Waldhausen donn\'e.
}\end{Resume}
%%%%%%%%%%%%%%%%%%%%%%%%%%%%%%%%%%%%%%%%%%%%%%%%%%%%%%%%%%%%
%%%  Le titre en Anglais  %%%
%%%%%%%%%%%%%%%%%%%%%%%%%%%%%
\selectlanguage{english}
\begin{Etitle}{%
The Novikov homology of graph links
}\end{Etitle}
%%%%%%%%%%%%%%%%%%%%%%%%%%%%%%%%%%%%%%%%%%%%%%%%%%%%%%%%%%%%
%%%  Abstract  %%%
%%%%%%%%%%%%%%%%%%
\begin{Abstract}{%
A graph multilink is a link with multiplicities in a homology 3-sphere whose 
exterior is a graph manifold. In this Note, we compute the Novikov homology of 
graph multilinks. As a corollary, we give a majoration for the number of Novikov modules
on a given graph link.
}\end{Abstract}
%%%%%%%%%%%%%%%%%%%%%%%%%%%%%%%%%%%%%%%%%%%%%%%%%%%%%%%%%%%%
%%%  Abridged English version  %%%
%%%%%%%%%%%%%%%%%%%%%%%%%%%%%%%%%%
%\AEv
%La version anglaise abr\'eg\'ee.(pas oblige)
%\par\medskip\centerline{\rule{2cm}{0.2mm}}\medskip
%\setcounter{section}{0}
%\selectlanguage{francais}
%%%%%%%%%%%%%%%%%%%%%%%%%%%%%%%%%%%%%%%%%%%%%%%%%%%%%%%%%%%%
%%%  Texte principal (en Francais)  %%%
%%%%%%%%%%%%%%%%%%%%%%%%%%%%%%%%%%%%%%%

\section{Introduction}
Soit $L=L_1\cup\dots \cup L_n$ un entrelacs orient\'e dans une 3-sph\`ere d'homologie $\Sigma$. Consid\'erons son ext\'erieur
$X=\Sigma-\text{int}\,{\cal N}(L)$, o\`u ${\cal N}(L)$ est un voisinage tubulaire de $L$. Les rev\^etements infini cycliques de
$X$ sont classifi\'es par
$$
[X,S^1]\simeq\Hom(\pi_1X,\ZZ)=\Hom(H_1X,\ZZ)\simeq H^1X\simeq H_1L=\bigoplus_{i=1}^n\ZZ[L_i].
$$
Ainsi, chaque rev\^etement infini cyclique de $X$ correspond \`a des multiplicit\'es
$\underline{m}=(m_1,\dots,m_n)$, o\`u $m_i$ est un entier associ\'e \`a $L_i$. Un entrelacs orient\'e $L$
muni d'une telle suite de multiplicit\'es $\underline{m}$ est appel\'e un {\it multi-entrelacs}, et not\'e
$L(\underline{m})=m_1L_1\cup\dots\cup m_nL_n$.

Si $\widetilde{X}(\mm)\stackrel{p}{\to}X$ d\'esigne le rev\^etement d\'etermin\'e par $\underline{m}$, le groupe ab\'elien
$H_1(\widetilde{X}(\mm))$ est muni d'une structure naturelle de $\ZZ[t,t^{-1}]$-module: on parle du {\it module d'Alexander
de\/} $L(\underline{m})$. Comme il existe toujours une matrice de pr\'esentation carr\'ee de ce module, le {\it polyn\^ome d'Alexander
de\/} $L(\underline{m})$ peut \^etre d\'efini comme le d\'eterminant d'une telle matrice. 

{\it L'homologie de Novikov de\/} $L(\underline{m})$ est le $\ZZ [\![t]\!][t^{-1}]$-module
$\widehat{H}_L(\mm) = H_1(\widetilde{X}(\mm))\otimes_{\ZZ [t,t^{-1}]}\ZZ [\![t]\!][t^{-1}]$.
Cette homologie, introduite par S.P.\,Novikov en 1981, lui permit de construire un analogue de la th\'eorie de Morse
pour les $1$-formes ferm\'ees non-exactes. 
\begin{exemple}
Un multi-entrelacs $L(\underline{m})$ est dit {\it fibr\'e\/} s'il existe une fibration localement triviale $X\stackrel{\phi}{\to}S^1$
dans la classe d'homotopie $\underline{m}$. Le polyn\^ome d'Alexander d'un tel multi-entrelacs est une unit\'e de l'anneau
$\ZZ [\![t]\!][t^{-1}]$; l'homologie de Novikov est donc triviale.
\end{exemple}

Voici une variation d'un probl\`eme pos\'e par Novikov dans [4]: \'etant donn\'e un entrelacs $L=L_1\cup\dots\cup L_n$ dans une 3-sph\`ere
d'homologie, comment se comporte l'homologie de Novikov $\widehat{H}_L(\mm)$ en fonction de $\mm\in\ZZ^n$? Nous donnons une
r\'eponse compl\`ete pour les entrelacs de Waldhausen, qu'il s'agit \`a pr\'esent de d\'efinir.
 
\section{Entrelacs de Waldhausen}
Soient $L'=L_0'\cup L_1'\cup \dots\cup L_n'$ et $L''=L_0''\cup L_1''\cup \dots\cup L_r''\/$ deux entrelacs d'ext\'erieurs $X'$
et $X''$ dans des sph\`eres d'homologie $\Sigma'$ et $\Sigma''$.
Choisissons des voisinages tubulaires ${\cal N}(L_0')$ et ${\cal N}(L_0'')$ munis de parall\`eles
et de m\'eridiens standards $P',M'\subset\partial{\cal N}(L_0')$, et $P'',M''\subset\partial{\cal N}(L_0'')$; posons
$\Sigma = (\Sigma'-{\cal N}(L_0'))\cup_{h}(\Sigma''-{\cal N}(L_0''))$, 
o\`u $h\colon \partial{\cal N}(L_0') \to \partial{\cal N}(L_0'')$ est un hom\'eomorphisme envoyant $P'$ sur $M''$ et $M'$ 
sur $P''$. Il est facile de v\'erifier que $\Sigma$ est une sph\`ere d'homologie. L'entrelacs $L=L_1'\cup \dots\cup L_n'\cup L_1''
\cup \dots\cup L_r''$ dans $\Sigma$ est appel\'e {\it l'\'epissure de $L'$ et $L''$ le long de $L_0'$ et $L_0''$\/}; il est not\'e
$$
L=L'\frac{}{\scriptstyle{L_0'}\quad \scriptstyle{L_0''}}L''.
$$
L'\'epissure $L'(\underline{m}')\frac{}{\scriptstyle{m_0'L_0'}\quad\scriptstyle{m_0''L_0''}}L''(\underline{m}'')$ de deux
multi-entrelacs est d\'efinie si et seulement si les multiplicit\'es $\underline{m}'\in \Hom(H_1X',\ZZ)$ et
$\underline{m}''\in \Hom(H_1X'',\ZZ)$ co\"\i ncident sur l'homologie du tore de recollement
$\partial{\cal N}(L_0')=\partial{\cal N}(L_0'')\subset\Sigma$. On aboutit aux \'equations
\begin{equation}\tag{$\ast$}
m_0' = \sum_{j=1}^r m_j''\ell k(L_0'',L_j'')\quad\text{ et }\quad m_0'' = \sum_{i=1}^n m_i'\ell k(L_0',L_i').
\end{equation}
Un {\it entrelacs simple\/} est un entrelacs irr\'eductible $L$ tel que tout tore incompressible dans l'ext\'erieur de $L$ est
parall\`ele au bord. Un {\it entrelacs seifertique\/} est un entrelacs dont l'ext\'erieur admet un feuilletage de Seifert.
Par le th\'eor\`eme de d\'ecomposition de Jaco-Shalen (voir [3]), tout entrelacs irr\'eductible $L$ dans une
3-sph\`ere d'homologie peut \^etre exprim\'e comme l'\'epissure d'un nombre fini d'entrelacs simples ou seifertiques, et
il existe une unique fa\c con minimale de le faire. On parle de la {\it d\'ecomposition de Jaco-Shalen de\/} $L$, les composantes
de cette d\'ecomposition \'etant les {\it \'el\'ements d'\'epissure de\/} $L$. Un {\it entrelacs de Waldhausen\/} est un entrelacs
dont tous les \'el\'ements d'\'epissure sont des entrelacs seifertiques.
\begin{exemple}
Les entrelacs de Waldhausen dans $S^3$ sont exactement les entrelacs toriques it\'er\'es, c'est-\`a-dire, les entrelacs
constructibles par cablage \`a partir d'entrelacs triviaux.
\end{exemple}

\begin{figure}
\begin{center}
\epsfig{figure=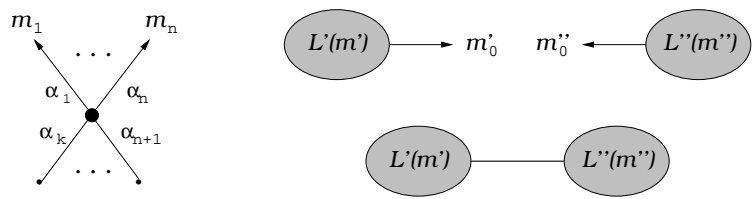,height=2.5cm}
\caption{Diagrammes d'\'epissure/{\it Splice diagrams}}
\label{fig:splice}
\end{center}
\end{figure}
\begin{figure}
\begin{center}
\epsfig{figure=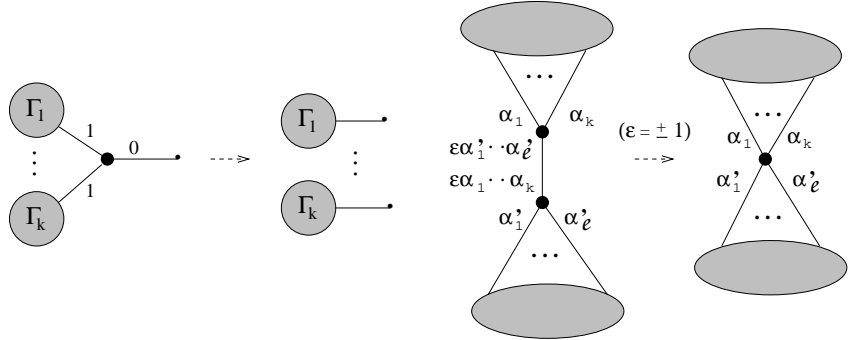,height=4.5cm}
\caption{Reductions de diagrammes/{\it Reduction of splice diagrams}}
\label{fig:red}
\end{center}
\end{figure}

Eisenbud et Neumann ont donn\'e dans [2] une classification des multi-entrelacs de Waldhausen au moyen de {\it diagrammes d'\'epissure},
d\'efinis comme suit. Tout multi-entrelacs seifertique peut \^etre repr\'esent\'e par un diagramme de la forme du premier dessin de la
Figure \ref{fig:splice}. Les $\alpha_i$ sont les {\it poids} associ\'es aux feuilles exceptionnelles du feuilletage de Seifert,
les pointes de fl\`eche symbolisent les composantes de l'entrelacs, et le sommet repr\'esente une feuille g\'en\'erique du feuilletage
de Seifert. L'\'epissure $L'(\underline{m}')\frac{}{\scriptstyle{m_0'L_0'}\quad\scriptstyle{m_0''L_0''}}L''(\underline{m}'')$
de deux multi-entrelacs est symbolis\'ee par les
diagrammes de $L'(\underline{m}')$ et de $L''(\underline{m}'')$ reli\'es le long des ar\^etes correspondant \`a $L_0'$ et \`a $L_0''$
(voir Figure \ref{fig:splice}). L'op\'eration inverse sur un diagramme consiste \`a couper une ar\^ete reliant deux sommets et \`a
ajouter deux pointes de fl\`eches munies des multiplicit\'es donn\'ees en ($\ast$); on parle de
{\it d\'ecomposition du diagramme\/}. 
Enfin, l'union disjointe de diagrammes repr\'esente la somme disjointe des multi-entrelacs.
Bien entendu, deux diagrammes d'\'epissure diff\'erents peuvent repr\'esenter le m\^eme multi-entrelacs. Par exemple,
les r\'eductions illustr\'ees en Figure \ref{fig:red} ne changent pas le multi-entrelacs repr\'esent\'e. Si aucune de ces deux
r\'eductions ne peut \^etre effectu\'ee, on dira que le diagramme est {\it minimal}.

Les diagrammes d'\'epissure sont parfaitement adapt\'es au calcul des coefficients d'enlacement. Etant donn\'es deux sommets ou
pointes de fl\`eche $v$ et $w$, le coefficient d'enlacement $\ell(v,w)$ des courbes correspondantes (feuille g\'en\'erique de Seifert
ou composante de l'entrelacs) est donn\'e par la formule suivante: soit $\sigma(v,w)$ la g\'eod\'esique du diagramme reliant $v$ et
$w$ ($v$ et $w$ inclus); alors, $\ell(v,w)$ est le produit de tous les poids des ar\^etes adjacentes \`a $\sigma(v,w)$ mais pas sur
$\sigma(v,w)$.

De plus, il est tr\`es facile de lire sur son diagramme minimal si un multi-entrelacs est fibr\'e ou non. Un
multi-entrelacs est fibr\'e si et seulement s'il est irr\'eductible et chacun de ses \'el\'ements d'\'epissure est fibr\'e. Un
multi-entrelacs seifertique est fibr\'e si et seulement si son
coefficient d'enlacement avec une feuille g\'en\'erique est non-nul. %(sauf si $L$ est l'entrelacs de Hopf, qui fibre toujours)
Ainsi, le multi-entrelacs repr\'esent\'e par le premier diagramme (minimal) de la Figure \ref{fig:splice} est fibr\'e si et seulement
si
$$
\sum_{i=1}^n m_i\alpha_1\cdots\widehat\alpha_i\cdots\alpha_k\neq 0.
$$
Par cons\'equent, les entrelacs de Waldhausen irr\'eductibles sont g\'en\'eriquement fibr\'es. En particulier, leur homologie de
Novikov est g\'en\'eriquement triviale.

Nous sommes enfin en mesure d'\'enoncer la formule pour l'homologie de Novikov d'un entrelacs de Waldhausen.

\section{R\'esultats}
Soit $L=L_1\cup\dots\cup L_n$ un entrelacs de Waldhausen donn\'e par un diagramme d'\'epissure minimal $\Gamma$, et soit
$\mm\in\ZZ^n-\{0\}$.
Il est facile de v\'erifier que les modules d'Alexander de $(L'(\underline{m}')\cup 0L_0')\frac{}{\quad}0L_0''$
et de $L'(\underline{m}')$ sont \'egaux; il en va donc de m\^eme pour les modules de Novikov, et nous allons supposer
que ces simplifications ont \'et\'e faites dans $\Gamma(\mm)$.

Chaque sommet de $\Gamma(\mm)$ correspond \`a un multi-entrelacs seifertique, fibr\'e ou non-fibr\'e; on parlera de sommet fibr\'e ou
non-fibr\'e. Par ailleurs, on dira qu'un sommet est un 0-sommet si toutes les multiplicit\'es du multi-entrelacs seifertique
correspondant sont nulles. Soit $\Gamma'(\mm')$ le sous-diagramme de $\Gamma(\mm)$ obtenu de la fa\c con suivante: on d\'ecompose
$\Gamma(\mm)$ le long des ar\^etes reliant un sommet fibr\'e et un sommet non-fibr\'e, et on
efface toutes les composantes connexes correspondant \`a des multi-entrelacs fibr\'es. Soient $c$ le nombre de composantes connexes
de $\Gamma$, $r$ le nombre de composantes connexes de $\Gamma'(\mm')$, $n$ le nombre de fl\`eches de $\Gamma'(\mm')$ et $k$ le nombre
de ses sommets. Pour chaque sommet $v$ de $\Gamma'(\mm')$, notons encore $\alpha(v)$ le produit des poids des fl\`eches de
$\Gamma'(\mm')$ adjacentes \`a $v$.
\begin{theoreme}
L'homologie de Novikov du multi-entrelacs $L(\mm)$ est la somme directe
\begin{itemize}
\item[-] d'un facteur libre de rang $n+c-r-k-1$, et 
\item[-] d'un facteur pr\'esent\'e par la $(k\times k)$-matrice $P=(p_{vw})$, avec $v$ et $w$ parcourant l'ensemble des sommets
non-fibr\'es, o\`u $p_{vw}=\frac{\ell(v,w)}{\alpha(w)}$ si
$v$ et $w$ appartiennent \`a la m\^eme composante connexe de $\Gamma'(\mm')$ et $v$ n'est pas un 0-sommet, et $p_{vw}=0$ sinon.
\end{itemize}
\end{theoreme}
Donnons un tr\`es bref aper\c cu des techniques utilis\'ees. La premi\`ere \'etape consiste \`a adapter les notions classiques de
surfaces et de formes de Seifert aux multi-entrelacs. Ces formes de Seifert ``g\'en\'eralis\'ees'' permettent de calculer
le module d'Alexander des multi-entrelacs seifertiques non-fibr\'es.
%Avec les notations du th\'eor\`eme, ce $\ZZ[t,t^{-1}]$-module est la somme directe d'un facteur
%libre de rang $n-k-1$ et d'un facteur pr\'esent\'e par la $(k\times k)$-matrice 
%$\bigl(\frac{\ell(v,w)}{\alpha(w)}{\scriptstyle (t^{d_v}-1)}\bigr)_{v,w}$, o\`u $d_v$ est le plus grand commun diviseur des 
%multiplicit\'es du multi-entrelacs correspondant \`a $v$ (voir [1] pour les d\'etails).
Ensuite, il s'agit de comprendre l'effet de l'\'epissure sur les modules consid\'er\'es. En travaillant sur $\ZZ[t,t^{-1}]$ (i.e: en
\'etudiant le module d'Alexander), il est extr\^emement difficile de donner une formule close. En revanche, le passage \`a
l'homologie de Novikov (en d'autres termes, la tensorisation par $\ZZ [\![t]\!][t^{-1}]$) permet de simplifier consid\'erablement
les calculs (voir [1] pour les d\'etails).
%L'\'epissure d'un multi-entrelacs de Waldhausen purement non-fibr\'e avec un multi-entrelacs fibr\'e ne change pas
%le module de Novikov. Comme $t^{d_v}-1$ est une unit\'e de l'anneau de Novikov \`a moins que $d_v=0$ (cas d'un 0-sommet),
%on obtient la formule du th\'eor\`eme.

\begin{corollaire}
Soit $L=L_1\cup\dots\cup L_n$ un entrelacs de Waldhausen. Il existe un nombre fini d'hyperplans dans $\ZZ^n$, d\'efinis par des
\'equations lin\'eaires homog\`enes \`a coefficients entiers, tels que $\widehat{H}_L(\mm)$ est constant sur chaque
strate de la ``stratification'' de $\ZZ^n$ par ces hyperplans.
\end{corollaire}

Voir [5], Th\'eor\`eme 3, p. 529 pour un r\'esultat analogue dans un autre contexte.

\begin{corollaire}
Soit un entrelacs de Waldhausen repr\'esent\'e par un diagramme d'\'epissure minimal \`a $k$ sommets et $c$ composantes connexes.
Alors, le nombre de modules de Novikov est born\'e par $3^k-2(k-c)$.
\end{corollaire}

Comme nous l'a fait remarquer A. Pazhitnov, ce ph\'enom\`ene de finitude du nombre de modules semble exceptionnel en homologie de Novikov.

\begin{exemple}
Notons $\hat\Lambda$ l'anneau de Novikov $\ZZ [\![t]\!][t^{-1}]$. L'entrelacs seifertique $L$ correspondant au
diagramme minimal donn\'e en Figure \ref{fig:splice} satisfait
$$
\widehat{H}_L(\mm)=
\begin{cases}
	0& \text{si $\mm\in\ZZ^n-V$},\\
	\hat\Lambda/(\alpha_{n+1}\cdots\alpha_k)\,\oplus\,{\hat\Lambda}^{n-2}& \text{si $\mm\in V-\{0\}$},\\				
	{\hat\Lambda}^n& \text{si $\mm=0$},
\end{cases}
$$
o\`u $V=\{(m_1,\dots,m_n)\in\ZZ^n\;|\;\sum_{i=1}^n m_i\alpha_1\cdots\widehat\alpha_i\cdots\alpha_k=0\}$.
\end{exemple}
%%%%%%%%%%%%%%%%%%%%%%%%%%%%%%%%%%%%%%%%%%%%%%%%%%%%%%%%%%%%
%%%  Remerciements  %%%
%%%%%%%%%%%%%%%%%%%%%%%
\Remerciements{Je souhaite remercier vivement Claude Weber, Andrei Pazhitnov et Mathieu Baillif
pour leur aide. C'est \'egalement un plaisir de t\'emoigner toute ma gratitude \`a Jerry Levine pour son accueil
et sa disponibilit\'e \`a l'Universit\'e de Brandeis.}
%%%%%%%%%%%%%%%%%%%%%%%%%%%%%%%%%%%%%%%%%%%%%%%%%%%%%%%%%%%%
%%%  Bibliographie %%%
%%%%%%%%%%%%%%%%%%%%%%

%
\end{document}